\documentclass[10pt]{article}
\usepackage{amstext}
\usepackage{amsfonts}
\usepackage{amssymb}
\usepackage{amsbsy}
\usepackage{latexsym}
\usepackage{xy}
\usepackage{hhline}
\xyoption{all}
\newcounter{num}[section] %

\newenvironment{theo}
{\refstepcounter{num}%
\bigskip\noindent{\bf Theorem~\arabic{section}.\arabic{num}. }\it}
\newenvironment{cor}
{\refstepcounter{num}%
\bigskip\noindent{\bf Corollary~\arabic{section}.\arabic{num}. }\it}
\newenvironment{lemma}
{\refstepcounter{num}%
\bigskip\noindent{\bf Lemma~\arabic{section}.\arabic{num}. }\it}
\newcommand{\example}
{\refstepcounter{num}%
\bigskip\noindent{\bf Example~\arabic{section}.\arabic{num}.}}

\newcommand{\remark}
{\refstepcounter{num}%
\bigskip\noindent{\bf Remark~\arabic{section}.\arabic{num}.}}

\newenvironment{proof}{\medskip\noindent{\it Proof. }}
{$\Box$ \bigskip}
\newenvironment{proof_of}[1]{\medskip\noindent{\it Proof #1}}
{$\Box$ \bigskip}

\newenvironment{eq}{\begin{equation}}{\end{equation}}

\newcommand{\Ref}[1]{(\ref{#1})}

\newcommand{\si}{\sigma}
\newcommand{\al}{\alpha}
\newcommand{\be}{\beta}

\newcommand{\la}{\lambda}
\newcommand{\de}{\delta}

\newcommand{\ov}[1]{\overline{#1}}
\newcommand{\un}[1]{{\underline{#1}} }

\newcommand{\tr}{\mathop{\rm tr}}

\newcommand{\mdeg}{\mathop{\rm mdeg}}
\newcommand{\hterm}{\mathop{\rm ht}}   % the highest term 
\newcommand{\Skew}{\mathop{\rm skew}}   

\newcommand{\Char}{\mathop{\rm char}}

\newcommand{\trdeg}{\mathop{\rm{tr.deg }}}
\newcommand{\algA}{\mathcal{A}}    %algebra A
\newcommand{\W}{\mathcal{W}} % the set of monomials in the generic matrices 

\newcommand{\FF}{{\mathbb{F}}}   % base field
\newcommand{\CC}{{\mathbb{C}}}   % base field
\newcommand{\NN}{{\mathbb{N}}}
\newcommand{\ii}{{\mathbb{I}}}   % (-1)^{1/2}

\newcommand{\Sp}{S\!p}

\begin{document}
\title{Orthogonal invariants of skew-symmetric matrices}
 \author{
A.A. Lopatin\thanks{Supported by DFG} \\
{\small\it Omsk Institute of Mathematics, SB RAS, Pevtsova street, 13, Omsk 644099, Russia} \\
{\small\it artem\underline{ }lopatin@yahoo.com} \\
%http://www.iitam.omsk.net.ru/\~{}lopatin/\\
}
\date{} % !!!
\maketitle

\begin{abstract} 
The algebra of invariants of $d$-tuples of $n\times n$ skew-symmetric matrices under the action of the orthogonal group by simultaneous conjugation is considered over an infinite field of characteristic different from two.  For $n=3$ and $d>0$ a minimal set of generators is established. A homogeneous system of parameters (i.e., an algebraically independent set such that the algebra of invariants is a finitely generated free module over subalgebra generated by this set) is described for $n=3$ and $d>0$, for $n=4$ and $d=2,3$, for $n=5$ and $d=2$.  
\end{abstract}



%SSSSSSSSSSSSSSSSSSSSSSSSSSSSSSSSSSSSSSSSSSSSSSSSSSSSSSSSSSSSSSSSSSSSSSSSSS
%SSSSSSSSSSSSSSSSSSSSSSSSSSSSSSSSSSSSSSSSSSSSSSSSSSSSSSSSSSSSSSSSSSSSSSSSSS
\section{Introduction}\label{section_intro}

We assume that $\FF$ is an infinite field of characteristic $\Char{\FF}$ different from two unless otherwise stated. All vector spaces, algebras, and modules are over $\FF$ and all algebras are associative. 

For $n>1$ and $d>0$ the orthogonal group $O(n)$ acts on $d$-tuples 
$$V=\FF^{n\times n}\oplus\cdots\oplus \FF^{n\times n}\text{ and }V_{-}=S_{-}\oplus\cdots \oplus S_{-}$$
of $n\times n$ matrices ($n\times n$ skew-symmetric matrices, respectively) over $\FF$ by the
diagonal conjugation, i.e.,
\begin{eq}\label{eq_diag_conj}
g\cdot (A_1,\ldots,A_d)=(g A_1 g^{-1},\ldots,g A_d g^{-1}),
\end{eq}%
where $g\in O(n)$ and $A_1,\ldots,A_d$ lay in $\FF^{n\times n}$ ($S_{-}$, respectively). The coordinate rings of the affine varieties $V$ and $V_{-}$ are 
the following algebras  %
$$R=\FF[V]=\FF[x_{ij}(k)\,|\,1\leq i,j\leq n,\, 1\leq k\leq d] \text{ and }
R_{-}=\FF[x_{ij}(k)\,|\,1\leq i<j\leq n,\, 1\leq k\leq d].$$ %
Denote by
$$X_k=\left(\begin{array}{ccc}
x_{11}(k) & \cdots & x_{1n}(k)\\
\vdots & & \vdots \\
x_{n1}(k) & \cdots & x_{nn}(k)\\
\end{array}
\right)
$$%
the $k^{\rm th}\!\!$ {\it generic} matrix and denote by $Y_{k}$ the $n\times n$ {\it skew-symmetric generic} matrix, i.e., 
$$(i,j)^{\rm th}\text{ entry of }Y_k =
\left\{\begin{array}{rl}
x_{ij}(k),& \text{if }  i<j\\
-x_{ji}(k),& \text{if }  i>j\\
0,& \text{otherwise}\\
\end{array}
\right..$$

The action of $O(n)$ on $V$ induces the action on $R$ as follows: $(g\cdot
f)(h)=f(g^{-1}\cdot h)$ for all $g\in O(n)$, $f\in \FF[V]$,
$h\in V$. In other words,  %
$$g\cdot x_{ij}(k)= (i,j)^{\rm th}\text{ entry of }g^{-1}X_k g.$$%
The algebra of {\it $O(n)$-invariants of matrices} is
$$R^{O(n)}=\{f\in \FF[V]\,|\,g\cdot f=f\;{\rm for\; all}\;g\in O(n)\}.$$
In the same way $O(n)$ acts on $R_{-}$ and we obtain $R_{-}^{O(n)}$, the algebra of {\it $O(n)$-invariants of skew-symmetric matrices}. 

Denote coefficients in the characteristic polynomial
of an $n\times n$ matrix $X$ by $\si_t(X)$, i.e., %
\begin{eq}\label{eq1_intro} 
\det(X+\la E)=\sum_{t=0}^{n} \la^{n-t}\si_t(X).
\end{eq}%
So, $\si_0(X)=1$, $\si_1(X)=\tr(X)$, and $\si_n(X)=\det(X)$. Part~a) of the following theorem was proven in~\cite{Zubkov99} and part~b) in~\cite{Lopatin_so_inv} (see also~\cite{ZubkovI}). 

\begin{theo}\label{theo_matrix}
Assume that $\Char{\FF}\neq2$ and $P$ is $R$ or $R_{-}$. Then the algebra of invariants  
$P^{O(n)}$ is generated by $\si_t(B)$ ($1\leq t\leq n$), where $B$ ranges over all monomials in 
\begin{enumerate}
\item[a)] $X_1,\ldots,X_d$, $X_1^T,\ldots,X_d^T$, if $P=R$; 

\item[b)] $Y_1,\ldots,Y_d$, if $P=R_{-}$. 
\end{enumerate}
Moreover, in both cases we can assume that $B$ is primitive, i.e., is not equal to a power of a shorter monomial. 
\end{theo}

\example\label{ex1} For $g\in O(n)$ we have $g\cdot\tr(Y_1 Y_2) = \tr(g^{-1} Y_1 g g^{-1} Y_2 g) = \tr(Y_1 Y_2)$. Therefore, $\tr(Y_1 Y_2)\in R_{-}^{O(n)}$.

\remark\label{remark_tr} In the case of a characteristic zero field it is enough to take
traces instead of $\si_t$, $1\leq t\leq n$, in the formulation of
Theorem~\ref{theo_matrix}.

\remark\label{remark_G} If $G$ is a classical linear group, i.e., $G$ belongs to the list $GL(n)$, $O(n)$, $\Sp(n)$, $SO(n)$, $SL(n)$, then we can define the algebra of invariants $R^{G}$ in the same way as above.  A generating set for the algebra $R^G$ is known for an arbitrary characteristic of $\FF$. It was established  in~\cite{Sibirskii_1968},~\cite{Procesi76},~\cite{Aslaksen95} for characteristic zero case and in~\cite{Donkin92a},~\cite{Zubkov99},~\cite{Lopatin_so_inv} for the general case. Relations between generators for $R^{GL(n)}$ were described in~\cite{Zubkov96} and for $R^{O(n)}$ in~\cite{Lopatin_Orel} (modulo free relations). Note that we always assume that  $\Char{\FF}\neq2$ in the case of $O(n)$ and $SO(n)$. 
\bigskip

Given $f\in R$, denote by $\deg{f}$ its {\it degree} and by $\mdeg{f}$ its {\it multidegree}, i.e.,
$\mdeg{f}=(t_1,\ldots,t_d)$, where $t_k$ is the total degree of the polynomial $f$ in $x_{ij}(k)$, $1\leq i,j\leq n$, and $\deg{f}=t_1+\cdots+t_d$. Since $\deg{\si_t(Z_1\cdots Z_s)}=ts$, where $Z_k$ is a generic or a transpose generic matrix, the algebra $R^{O(n)}$ as well as $R_{-}^{O(n)}$ has $\NN$-grading by degrees and $\NN^d$-grading by multidegrees, where $\NN$ stands for non-negative integers.

Assume that $P$ is $R$ or $R_{-}$. By the Noether normalization lemma,  
$P^{O(n)}$ contains a {\it homogeneous} (with respect to $\NN$-grading)
{\it system of parameters} (shortly h.s.p.), i.e., a set $f_1,\ldots,f_s$ of
algebraically independent elements such that $P^{O(n)}$ is a finitely generated $\FF[f_1,\ldots,f_d]$-module. Moreover, since $P^{O(n)}$ is a Cohen-Macaulay algebra (see~\cite{Hashimoto_2001}),  $P^{O(n)}$ is a free $\FF[f_1,\ldots,f_s]$-module for any h.s.p.~$f_1,\ldots,f_s$. A system of
parameters for $R^{GL(n)}$ was constructed for $n=2$ and any $d$
(see~\cite{Teranishi_1988} and~\cite{DKZ_2002}), $n=3,4$ and $d=2$
(see~\cite{Teranishi_1986}), and $n=d=3$ (see~\cite{Lopatin_Sib}). In this paper we proved the following result.

\begin{theo}\label{theo_hsp}
The sets 
\begin{enumerate} 
\item[a)] $\si_2(Y_i),\, 1\leq i\leq d$; $h_r=\sum \tr(Y_i Y_j),\, 3\leq r\leq 2d+1$, 

where $n=3$, $d>0$, and the sum in the definition of $h_r$ ranges over all $1\leq i<j\leq d$ with $i+j=r$;

\item[b)] $\si_2(Y_1)$, $\si_2(Y_2)$, $\det(Y_1)$, $\det(Y_2)$, $\tr(Y_1 Y_2)$, $\tr(Y_1^2 Y_2^2)$, 

where $n=4$ and $d=2$; 

\item[c)] $\si_2(Y_i)$, $\det(Y_i)$, $1\leq i\leq 3$;  $\tr(Y_i Y_j)$, $\tr(Y_i^2 Y_j^2)$, $1\leq i<j\leq 3$, 

where $n=4$ and $d=3$; 

\item[d)] $\si_2(Y_1)$, $\si_2(Y_2)$, $\si_4(Y_1)$, $\si_4(Y_2)$,  $\tr(Y_1 Y_2)$,  $\tr(Y_1^2 Y_2^2)$, $\tr(Y_1^3 Y_2)$,  $\tr(Y_1 Y_2^3)$, $\tr(Y_1^4 Y_2^2)$,  $\tr(Y_1^2 Y_2^4)$,   

where $n=5$ and $d=2$
\end{enumerate}
are homogeneous systems of parameters for $R_{-}^{O(n)}$ for the corresponding $n$ and $d$.
\end{theo}
\bigskip 

A {\it minimal homogeneous set of generators} (shortly m.h.s.g.) for $P^{O(n)}$ is a minimal (by inclusion) $\NN^d$-homogeneous set generating the algebra $P^{O(n)}$ over $\FF$.  A
m.h.s.g.~for $R^{GL(n)}$ is known for $n=2$ (see~\cite{Sibirskii_1968},~\cite{Procesi_1984},~\cite{DKZ_2002}) and $n=3$
(see~\cite{Lopatin_Comm1},~\cite{Lopatin_Comm2}). In characteristic zero case a
m.h.s.g.~for $R^{GL(n)}$ was also established for $n=4$ and $d=2$ (see~\cite{Drensky_Sadikova_4x4}). In this paper we obtained the following result.

\begin{theo}\label{theo_mgs}
The set
$$\si_2(Y_i),\, 1\leq i\leq d;\;\; \tr(Y_i Y_j), 1\leq i<j\leq d;\;\; \tr(Y_i Y_j Y_k),\, 1\leq i<j<k\leq d$$
is a minimal homogeneous set of generators for $R_{-}^{O(3)}$ for all $d>0$.
\end{theo}
\bigskip

The paper is organized as follows. In Section~\ref{section_notations} we introduce some notations and formulate key Hilbert Theorem. Then we consider the case of several $2\times 2$ matrices and the case of one $n\times n$ matrix. 

Proofs of Theorems~\ref{theo_hsp} and~\ref{theo_mgs} are given in Sections~\ref{section_n3},~\ref{section_n4},~\ref{section_n5} for $n=3,4,5$, respectively. It is not difficult to see that if Theorem~\ref{theo_hsp} holds over the algebraic closure of $\FF$, then it holds over $\FF$. 
%...It is really easy! 
So during the proof of Theorem~\ref{theo_hsp} we can assume that $\FF$ is algebraically closed.

\section{Notations and preliminaries}\label{section_notations}

Given an $\NN^d$-graded algebra $\algA$, denote by $\algA^{\#}$ the subalgebra generated by homogeneous elements of positive degree. A set $\{a_i\} \subseteq \algA$ is a m.h.s.g.~if and only if the $a_i$'s are $\NN^d$-homogeneous and $\{\ov{a}_i\}$ is a basis for  $\ov{\algA}={\algA}/{(\algA^{\#})^2}$. If we consider $a\in\algA$ as an element of $\ov{\algA}$, then we usually omit the bar and write $a\in\ov{\algA}$ instead of $\ov{a}$. An element $a\in \algA$ is called {\it
decomposable} if $a=0$ in $\ov{\algA}$. In other words, a decomposable
element is equal to a polynomial in elements of strictly lower degree.

Denote by $\W^{\#}$ ($\W_{-}^{\#}$, respectively) the monoid without unity generated by generic matrices $X_1,\ldots,X_d$, $X_1^T,\ldots,X_d^T$ ($Y_1,\ldots,Y_d$, respectively).  
We write $\W$ for  $\W^{\#}\sqcup\{E\}$, i.e., we endow $\W$ with the unity. Similarly, we denote $\W_{-} =\W_{-}^{\#}\sqcup\{E\}$. Given $f\in R_{-}^{O(n)}$ and $\un{A}\in V_{-}$, we write $f(\un{A})$ for the image of $\un{A}$ with respect to $f:V_{-}\to \FF$. In other words, $f(\un{A})$ is the result of substitution $X_k\to A_k$ in $f$ ($1\leq k\leq d$)

Given $a_1,\ldots,a_s\in\FF$, where $s=n(n-1)/2$, we write $\Skew(a_1,\ldots,a_s)$ for the following skew-symmetric matrix: 
$$\left(\begin{array}{ccccc}
0 & a_1 & \cdots & a_{n-2} & a_{n-1}\\
-a_1 &0 & \cdots & a_{2n-4}& a_{2n-3}\\
\vdots &\vdots & & \vdots &\vdots \\
-a_{n-2} & -a_{2n-4} & \cdots& 0 & a_s \\
-a_{n-1} & -a_{2n-3} & \cdots& -a_s & 0 \\
\end{array}
\right)\in S_{-}.$$

Given $A,B\in S_{-}$, we write $A\sim B$ if and only if there is an $g\in O(n)$ with $g A g^{-1}=B$. If the field $\FF$ is algebraically closed, then  we denote by $\ii$ some element of $\FF$ satisfying $\ii^2=-1$.

The following theorem was proven by Hilbert~\cite{Hilbert_old} in characteristic zero case  (for the contemporary reprint see~\cite{Hilbert_new}). Its proof from~\cite{Kraft} is valid for an arbitrary $\Char{\FF}$.

\begin{theo}\label{theo_Hilbert}
Let an algebraic group $G$ act regularly on some affine variety $W$.
This action induces the action of $G$ on the coordinate algebra $\FF[W]$ that
consists of polynomial maps from $W$ into $\FF$. Let invariants
$f_1,\ldots,f_s\in \FF[W]^{G}$ have the following property:
if $f_1(w)=\cdots=f_s(w)=0$, where $w\in W$, then for each 
invariant $f\in (\FF[W]^G)^{\#}$ we have $f(w)=0$. Then the
algebra of invariants $\FF[W]^G$ is a finitely generated module over its subalgebra generated by $f_1,\ldots,f_s$.
\end{theo}
\bigskip

Denote by  $\trdeg{R_{-}^{O(n)}}$ the {\it transcendence degree} of $R_{-}^{O(n)}$, i.e., the cardinality of its h.s.p.
We will use the following lemma to construct h.s.p.-s.

\begin{lemma}\label{lemma_key}
Let $f_1,\ldots,f_s\in R_{-}^{O(n)}$ be a set of $\NN$-homogeneous elements such that
\begin{enumerate}
\item[a)] if $f_1(\un{A})=\cdots=f_s(\un{A})=0$ for an $\un{A}=(A_1,\ldots,A_d)\in V_{-}$, then $f(\un{A})=0$ for all $f\in (R_{-}^{O(n)})^{\#}$; 

\item[b)] $s=n(n-1)(d-1)/2$.
\end{enumerate}
Then $f_1,\ldots,f_s$ is a h.s.p.~for $R_{-}^{O(n)}$.
\end{lemma}
\begin{proof}
It is well known (for example, see~\cite{Kraft}) that $\trdeg{R_{-}^{O(n)}}\geq \trdeg{R_{-}} - \dim{O(n)} = n(n-1)(d-1)/2=s$. On the other hand, Hilbert Theorem and condition~a) imply that $\trdeg{R_{-}^{O(n)}}\leq s$ and the required is proven. 
\end{proof}

Let us remark that $s$ from part b) of the previous lemma is not always equal to the cardinality of a h.s.p.~for $R_{-}^{O(n)}$. As examples, see below algebras of invariants from Lemmas~\ref{lemma_n2} and~\ref{lemma_d1}.

The next lemma follows from the fact that $\si_t(Y_1)=\si_t(Y_1^T)=\si_t(-Y_1)$.

\begin{lemma}\label{lemma_key2}
Let $n>1$. Then 
\begin{enumerate}
\item[a)] $\si_t(Y_1)=0$ for odd $t$ with $1\leq t\leq n$; 

\item[b)] $\tr(Y_1 Y_3 Y_2)=-\tr(Y_1 Y_2 Y_3)$; in particular, $\tr(Y_1^2 Y_2)=0$;

\item[c)] $2\si_2(U)=-\tr(U^2) + \tr(U)^2$ for all $U\in\W_{-}^{\#}$. 
\end{enumerate}
\end{lemma}

We conclude this section with consideration of two trivial cases.

\begin{lemma}\label{lemma_n2}
Let $n=2$ and $d>0$. Then 
\begin{enumerate}
\item[$\bullet$] $\si_2(Y_i), 1\leq i\leq d;\; \tr(Y_i Y_j), 1\leq i<j\leq d$ 

is a minimal homogeneous set of generators for $R^{O(2)}_{-}$;

\item[$\bullet$] $\si_2(Y_i), 1\leq i\leq d$   

is a homogeneous system of parameters for $R^{O(2)}_{-}$.
\end{enumerate} 
\end{lemma}
\begin{proof}
For short, we write $x_i$ for $x_{12}(i)\in R_{-}$. The required statement follows from the fact that $\si_2(Y_i)=x_i^2$ and 
$$\tr(Y_{i_1}\cdots Y_{i_s}) =
\left\{\begin{array}{rl}
2 x_{i_1}\cdots x_{i_s},& \text{if $s$ is even}  \\
0,& \text{otherwise}\\
\end{array}
\right.
$$
for $1\leq i,i_1,\ldots,i_s\leq d$.
\end{proof}

\begin{lemma}\label{lemma_d1}
Let $n>1$ and $d=1$. Then $R_{-}^{O(n)}\simeq \FF[\si_{2k}(Y_1)\,|\,1\leq 2k \leq n]$. 
\end{lemma}
\begin{proof}
For short, we write $x_{ij}$ for $x_{ij}(1)\in R_{-}$. Endow the set of monomials in $x_{ij}$, $1\leq i<j\leq n$, with the following partial lexicographical
order:
\begin{enumerate}
\item[$\bullet$] $x_{ij}> x_{pq}$ if $i<p$ or $i=p$ and $j<q$; 

\item[$\bullet$] $x_{i_1 j_1}\cdots x_{i_r j_r} > x_{p_1 q_1}\cdots x_{p_s q_s}$ if $x_{i_1 j_1} = x_{p_1 q_1},\ldots, x_{i_l j_l} = x_{p_l q_l}$ and $x_{i_{l+1} j_{l+1}} = x_{p_{l+1} q_{l+1}}$ for some $0\leq l<\min\{r,s\}$.
\end{enumerate}
Note that two words with different degrees can be incomparable. For an $\NN$-homogeneous $f\in R_{-}$ denote by $\hterm(f)$ the highest term of $f$.  In other words, if $f=\sum_i \al_i f_i$, where $\al_i\in\FF$, $\al_i\neq0$, and $f_i$ is a monomial, then $\hterm(h)=\max\{f_i\}$. Note that $\hterm(f)$ is well defined. 

For $1\leq 2k\leq n$ we have $\hterm(\si_{2k}(Y_1))=x_{12}^2\cdots x_{2k-1,2k}^2$. If elements $\si_{2k}(Y_1)$, $1\leq 2k\leq n$, are not algebraically independent over $\FF$, then their highest terms are also not algebraically independent; a contradiction. Theorem~\ref{theo_matrix} completes the proof. 
\end{proof}

\section{Canonical forms}\label{section_canonical}

In this section we assume that $\FF$ is algebraically closed. 
Given $p>0$, we consider the following $p\times p$ matrices:
$$A^{(p)}=\left(
\begin{array}{cccccc}
0     &1     &\cdot &\cdot &      &0     \\
-1    &0     &\cdot &      &      &\cdot \\
\cdot &\cdot &\cdot &\cdot &      &\cdot \\
\cdot &      &\cdot &\cdot &\cdot &\cdot \\
\cdot &      &      &\cdot &0     &1     \\
0     &\cdot &\cdot &\cdot &-1    &0     \\
\end{array}
\right),\quad
B^{(p)}=\left(
\begin{array}{cccccc}
0     &\cdot &\cdot &\cdot &1     &0     \\
\cdot &      &      &\cdot &0     &1     \\
\cdot &      &\cdot &\cdot &\cdot &\cdot \\
\cdot &\cdot &\cdot &\cdot &      &\cdot \\
1     &0     &\cdot &      &      &\cdot \\
0     &1     &\cdot &\cdot &\cdot &0     \\
\end{array}
\right),\text{ and }
$$
$$
C^{(p)}=\left(
\begin{array}{cccc}
0    &\cdots &0   &1     \\
0    &\cdots &1   &0     \\
     &\vdots &    &     \\
1    &\cdots &0   &0     \\
\end{array}
\right),
$$%
where $A^{(p)}$ and $B^{(p)}$ have exactly $2(p-1)$ non-zero elements. 
For $\la\in\FF$ we consider the following $p$ dimensional vectors:   
$$U_{\la}^{(p)}=\left(
\begin{array}{c}
\la    \\
0    \\
\vdots \\
0    \\
\end{array}
\right)\text{ and } 
V_{\la}^{(p)}=\left(
\begin{array}{c}
0    \\
\vdots \\
0    \\
\la \\
\end{array}
\right).
$$
For $\mu\in\FF$ we set  
$
K_{\mu}^{(2p)}=\frac{1}{2}\left(
\begin{array}{cc}
A^{(p)}    & \ii B^{(p)}+2\mu C^{(p)}\\
-\ii B^{(p)}-2\mu C^{(p)} & -A^{(p)} \\
\end{array}
\right)$ and 
$$K^{(2p+1)}=\frac{1}{2}\left(
\begin{array}{ccc}
A^{(p)}    & V^{(p)}_{1+\ii}& \ii B^{(p)}\\
 (V^{(p)}_{-1-\ii})^T&0&  (U^{(p)}_{-1+\ii})^T\\
-\ii B^{(p)}&  U^{(p)}_{1-\ii}& -A^{(p)} \\
\end{array}
\right).
$$%
We write $0^{(p)}$ for zero $p\times p$ matrix. The following result is well known in characteristic zero case (see~\cite{Wellstein_1930}). 

\begin{theo}\label{theo_canonical}
Assume that $A\in S_{-}$ is an $n\times n$ matrix. Then there is a $g\in O(n)$ such that $g A g^{-1}$ is a block-diagonal matrix, where non-zero blocks are $K^{(2p)}_{\mu}$, $K^{(2p+1)}$ for $p>0$, $\mu\in\FF$. 
\end{theo}
\bigskip

The proof of this theorem for $\FF=\CC$ can be found in~\cite{Gantmacher}. To adopt this proof for the general case it is enough to notice that the following lemma holds over an algebraically closed field of arbitrary characteristic. For example, this result follows from Theorem~2.6 of~\cite{DW_GenQuiv_2002}.  

\begin{lemma}\label{lemma_orbits}
Assume that $A,B\in S_{-}$ are $n\times n$ matrices and there exists a $g\in GL(n)$ such that $gAg^{-1}=B$. Then  $A\sim B$.
\end{lemma}
\bigskip

The next lemma is a corollary of Theorem~\ref{theo_canonical}.

\begin{lemma}\label{lemma_canon}
Assume that $A\in S_{-}$ is a non-zero $n\times n$ matrix and $\si_t(A)=0$ for all $1\leq t\leq n$. Then $A \sim B$, where 
\begin{enumerate}
\item[$\bullet$] if $n=3$, then $B=K^{(3)}=\frac{1}{2}\left(
\begin{array}{ccc}
0&1+\ii&0\\
-1-\ii&0&-1+\ii\\
0&1-\ii&0\\
\end{array}
\right);$

\item[$\bullet$] if $n=4$, then $B$ is $K^{(3)}\oplus 0^{(1)}$ or $K^{(4)}_0$;

\item[$\bullet$] if $n=5$, then $B$ one of the following matrices: $K^{(3)}\oplus 0^{(2)}$, $K^{(4)}_0\oplus 0^{(1)}$, $K^{(5)}$.
\end{enumerate}
\end{lemma}

\section{The case of $n=3$}\label{section_n3}

In this section we assume that $n=3$ and $d>0$.

\begin{lemma}\label{lemma_n3}
For $s\geq4$ and $1\leq i_1,\ldots,i_s\leq d$ we have $\tr(Y_{i_1}\cdots Y_{i_s})$ is decomposable in $R_{-}^{O(3)}$.
\end{lemma}
\begin{proof}
Let $s=4$. Since $4\tr(Y_1 Y_2 Y_3 Y_4) = \tr(Y_1 Y_2) \tr(Y_3 Y_4) + \tr(Y_1 Y_4) \tr(Y_2 Y_3)$, we have $\tr(Y_1 Y_2 Y_3 Y_4)\equiv 0$.

Let $s>4$. By formula~(5) and Theorem~5.1 from~\cite{Lopatin_O3}, $\tr((X_1-X_1^T)\cdots (X_4-X_4^T) U)$ is decomposable in $R^{O(3)}$, where $U\in\W^{\#}$. 
Thus, for $V\in\W_{-}^{\#}$ we have $2^4 \tr(Y_1 Y_2 Y_3 Y_4 V)\equiv 0$ in $R_{-}^{O(3)}$. The required is proven.   
\end{proof}

\begin{proof_of}{of Theorem~\ref{theo_mgs}.} Denote by $H$ the set from the formulation of the theorem.  Assume that $X,Y,Z$ are some skew-symmetric  generic matrices and $U\in\W_{-}^{\#}$. Using Theorem~\ref{theo_matrix} together with Lemmas~\ref{lemma_key2} and~\ref{lemma_n3} we obtain that $H$ generates $R_{-}^{O(3)}$. Since elements of $H$ belong to pairwise different $\NN^d$-homogeneous components of $R_{-}^{O(3)}$, to complete the proof it is enough to show that $\si_2(Y_1)$, $\tr(Y_1 Y_2)$, and $\tr(Y_1 Y_2 Y_3)$ are indecomposable.

Let $\si_2(Y_1)\equiv0$. Then $\si_2(Y_1)=\al \tr(Y_1)^2$ for an $\al\in\FF$. Therefore,  $\si_2(Y_1)=0$ (see part~a) of Lemma~\ref{lemma_key2}); a contradiction.  Similarly, we obtain that $\tr(Y_1 Y_2)\not\equiv0$.
 
Let $\tr(Y_1 Y_2 Y_3)\equiv0$. Then there are $\al,\be_1,\be_2,\be_3\in\FF$ such that $$\tr(Y_1 Y_2 Y_3)=\al \tr(Y_1) \tr(Y_2) \tr(Y_3) +
\be_1 \tr(Y_1) \tr(Y_2 Y_3) + \be_2 \tr(Y_2) \tr(Y_1 Y_3) + \be_3 \tr(Y_3) \tr(Y_1 Y_2).$$ Thus, $\tr(Y_1 Y_2 Y_3)=0$; a contradiction. The proof is completed. 
\end{proof_of}

In the rest of this section we assume that $\FF$ is an algebraically closed field.

\begin{lemma}\label{lemma_n3_trABC}
Let $A_1,A_2,A_3\in S_{-}$ and $\si_2(A_i)=0$ $(1\leq i\leq 3)$, $\tr(A_i A_j)=0$ ($1\leq i<j\leq 3$). Then $\tr(A_1 A_2 A_3)=0$. 
\end{lemma}
\begin{proof}
Denote $A_i=\Skew(a_i,b_i,c_i)$ for $1\leq i\leq3$.

Let $c_2\neq0$. Using $\tr(A_1 A_2)=0$ and $\tr(A_2 A_3)=0$ we obtain 
$$
c_1=-(a_1 a_2 + b_1 b_2)/c_2,\quad c_3=-(a_2 a_3 + b_2 b_3)/c_2,
$$%
respectively. Thus $\tr(A_1 A_2 A_3)=(a_3 b_1 - a_1 b_3) \si_2(A_2)/c_2 = 0$.

Let $c_2=0$. If $b_2\neq0$, then using $\tr(A_1 A_2)=\tr(A_2 A_3)=0$ we obtain $b_1 = -a_1 a_2/b_2$ and $b_3 = -a_2 a_3/b_2$; hence $\tr(A_1 A_2 A_3)=-(a_3 c_1 - a_1 c_3) \si_2(A_2)/b_2 = 0$.
If $b_2=0$, then $\si_2(A_2)=0$ implies $A_2=0$. 
\end{proof}

\begin{proof_of}{of part~a) of Theorem~\ref{theo_hsp}.} If $d=1$, then see Lemma~\ref{lemma_d1}.

Let $d>1$. Consider $A_1,\ldots,A_d\in S_{-}$ such that $\si_2(A_i)=0$ ($1\leq i\leq d$) and $h_r(\un{A})=0$ ($3\leq r\leq 2d+1$). Lemma~\ref{lemma_key} together with Theorem~\ref{theo_mgs} and Lemma~\ref{lemma_n3_trABC} implies that the required statement follows from 
\begin{eq}\label{eq1}
\tr(A_i A_j)=0 \text{ for all }1\leq i<j\leq d.
\end{eq} 
\indent{}We prove~\Ref{eq1} by induction on $2\leq l\leq d$, where we assume that $j\leq l$ in~\Ref{eq1}. If $l=2$, then $h_3(\un{A})=\tr(A_1 A_2)=0$.

Let $l\geq3$. By induction hypothesis, we have $\tr(A_i A_j)=0$ for all $1\leq i<j<l$. Let  $A_1=\cdots=A_{k-1}=0$ for $1\leq k< l$ and $A_k\neq0$. Note that if the mentioned $k$ does not exist, then $\tr(A_i A_l)=0$ for all $i<l$ and we obtain the required. 

We consider $i$ satisfying $k< i\leq l$. By Lemma~\ref{lemma_canon}, we can assume that $A_k=\Skew(1+\ii,0,-1+\ii)$ and $A_i=\Skew(a_i,b_i,c_i)$. By induction hypothesis and  $h_{k+l}(\un{A})=\tr(A_k A_l)=0$, we have $\tr(A_k A_i)=-2(1+\ii)a_i  -2(-1+\ii)c_i=0$. Therefore, $c_i=-\ii a_i$ and $\si_2(A_i)=0$ implies $b_i=0$. Thus, we have the equality $\tr(A_i A_l)=-2(a_i a_l+c_i c_l)=0$ which completes the proof. 
\end{proof_of}

\section{The case of $n=4$}\label{section_n4}

In this section we assume that $n=4$ and $\FF$ is algebraically closed. For every $d>0$ denote by $H_d$ the following set of invariants:
$$\si_2(Y_i),\, \det(Y_i),\, 1\leq i\leq d;\;\;  \tr(Y_i Y_j),\, \tr(Y_i^2 Y_j^2),\, 1\leq i<j\leq d.$$ 
We set $Q_1=K^{(3)}\oplus 0^{(1)}=\frac{1}{2}\Skew(1+\ii, 0,0, -1+\ii,0,0)$ and $Q_2=K^{(4)}_0=\frac{1}{2}\Skew(1,\ii,0,0,\ii,-1)$.

\begin{proof_of}{of part~b) of Theorem~\ref{theo_hsp}.} Let $d=2$. We assume that $A_1,A_2\in S_{-}$ satisfy 
\begin{eq}\label{eq2}
\si_2(A_i)=\det(A_i)=0\text{ for }i=1,2\text{ and }
\end{eq}\vspace{-5mm}%???...
\begin{eq}\label{eq3}
\tr(A_1 A_2)=\tr(A_1^2 A_2^2)=0.
\end{eq}% 
To prove the required statement, it is enough to show the following condition is valid (see Lemma~\ref{lemma_key}):
\begin{eq}\label{eq_cond}
f(A_1,A_2)=0\text{ for all }f\in (R_{-}^{O(4)})^{\#}.
\end{eq}% 
By Lemma~\ref{lemma_canon}, we can assume that $A_1\in\{0,Q_1,Q_2\}$ and $A_2=\Skew(a_2,b_2,c_2,d_2,e_2,f_2)$. 

If $A_1=0$, then condition~\Ref{eq_cond} follows from Lemma~\ref{lemma_d1}. 

Let $A_1=Q_1$. Equalities~\Ref{eq3} imply $a_2=-\ii\, d_2$, $c_2=\ii f_2$ and equalities~\Ref{eq2} imply $b_2=e_2=0$. Thus, $A_1^2 A_2 = A_2^2 A_1 = 0$ and $A_1 A_2 A_1 = A_2 A_1 A_2 = 0$. Lemma~\ref{lemma_key2} together with $\si_3(A_1 A_2)=0$ implies that condition~\Ref{eq_cond} holds.

Let $A_1=Q_2$. Equalities~\Ref{eq3} imply $a_2=-\ii\,(b_2 + e_2) + f_2$ and equalities~\Ref{eq2} imply $d_2=c_2$ and $c_2^2=(b_2 + \ii f_2)(e_2 + \ii f_2)$. Thus, $A_2^2 A_1 = A_1 A_2 A_1 = A_2 A_1 A_2 = 0$. Lemma~\ref{lemma_key2} together with $A_1^2=0$ and $\si_3(A_1 A_2)=0$ implies that condition~\Ref{eq_cond} holds.
\end{proof_of}

\begin{lemma}\label{lemma_n4_ABC}
Assume that for $A_1,A_2,A_3\in S_{-}$ we have $f(A_1,A_2,A_3)=0$ for all $f\in H_3$. Then $A_1 A_2 A_3 = 0$.  
\end{lemma}
\begin{proof} Obviously, without loss of generality we can assume that $A_i\neq0$ for $i=1,2,3$.

Let $A_i\sim Q_1$ for some $i$. Then without loss of generality we can assume that $A_1= Q_1$ and $A_i=\Skew(a_i,\ldots,f_i)$ for $i=2,3$. Similarly to the proof of part~b) of Theorem~\ref{theo_hsp}, we obtain $a_i=-\ii\, d_i$, $c_i=\ii f_i$, and $b_i=e_i=0$ ($i=2,3$). Thus, it is not difficult to see that $A_{\si(1)} A_{\si(2)} A_{\si(3)}=0$ for all $\si\in S_3$.

Let $A_i\sim Q_2$ for all $1\leq i\leq 3$. Consider $g\in O(4)$ satisfying $g A_1 g^{-1} = Q_2$. Denote $g A_i g^{-1}=B_i=\Skew(a_i,\ldots,f_i)$ for all $i$. In particular, $B_1=Q_2$. Similarly to the proof of part~b) of Theorem~\ref{theo_hsp}, we obtain $a_i=-\ii (b_i + e_i) + f_i$ and $d_i=c_i$ ($i=2,3$). Hence we have the equality $B_1 B_i = B_i B_1$ which implies $A_1 A_i = A_i A_1$. Similarly, we obtain $A_i A_j = A_j A_i$ for all $1\leq i<j\leq 3$. Thus, $B_i B_j = B_j B_i$.

For some $b_{11},b_{12},b_{22}\in\FF$ and $q=\tr(B_2 B_3)/2$ we have 
$$2 B_1 B_2 B_3 =\left(\begin{array}{cccc}
b_{11}               &   b_{12}  & \ii b_{12} & -\ii b_{11} \\
b_{12} - q           &  b_{22}   & \ii b_{22} & -\ii b_{12} + \ii q \\
\ii b_{12} - \ii  q  &\ii b_{22}&     -b_{22} & b_{12} - q \\
-\ii b_{11}          &-\ii b_{12}&     b_{12} & -b_{11} \\
\end{array}
\right).$$
Since $q=0$, we obtain $(B_1 B_2 B_3)^T = B_1 B_2 B_3$. On the other hand, $(B_1 B_2 B_3)^T=-B_3 B_2 B_1 = - B_1 B_2 B_3$. Thus, $B_1 B_2 B_3=0$ and we can see that $A_1 A_2 A_3 = 0$.

By Lemma~\ref{lemma_canon}, we have considered all possibilities for $A_1,A_2,A_3$. The proof is completed. 
\end{proof}

\begin{cor}\label{cor1}
For all $d>0$ the algebra $R_{-}^{O(4)}$ is integral over the subalgebra generated by $H_d$.
\end{cor}
\begin{proof}
Let $U,V\in\W^{\#}$. Since $\si_6(U+V)=0$ is a relation for $R^{O(4)}$, we obtain $\si_3(UV)\equiv -\tr(U^3 V^3) - \tr(U^2 V^2 U V) - \tr(U^2 V U V^2)$ in $R^{O(4)}$ (for details, see~\cite{Zubkov96} or~\cite{Lopatin_Orel}). Thus, the same equality holds in $R_{-}^{O(4)}$ for $U,V\in\W_{-}^{\#}$. Lemmas~\ref{lemma_key2} and~\ref{lemma_n4_ABC} together with Hilbert Theorem complete the proof.  
\end{proof}

The proof of part~c) of Theorem~\ref{theo_hsp} follows from Lemma~\ref{lemma_key} and Corollary~\ref{cor1}.

\section{The case of $n=5$}\label{section_n5}

In this section we assume that $n=5$ and $\FF$ is algebraically closed.  
We set %
$Q_1=K^{(3)}\oplus 0^{(2)} = \frac{1}{2}\Skew(1+\ii, 0,0, 0, -1+\ii,0,0,0,0,0)$,  %
$Q_2=K^{(4)}_0\oplus 0^{(1)}=\frac{1}{2}\Skew(1,\ii,0,0,0,\ii,0,-1,0,0)$, and    %
$Q_3=K^{(5)} =               \frac{1}{2}\Skew(1,0,\ii,0,1+\ii,0,\ii,-1+\ii,0,-1)$. %
Denote by $H$ the set from the formulation of part~d) of Theorem~\ref{theo_hsp}. 

\begin{lemma}\label{lemma_n5_1}
Assume that for $A_1,A_2\in S_{-}$ we have $f(A_1,A_2)=0$ for all $f\in H$ and  \begin{enumerate}
\item[a)] $\tr(A_1^3 A_2^3)=0$;

\item[b)] $\tr(A_1^{i_1} A_2^{j_1} \cdots A_1^{i_r} A_2^{j_r})=0$, where $r\geq2$ and $1\leq i_1,j_1,\ldots,i_r,j_r\leq 4$.
\end{enumerate} 
Then $f(A_1,A_2)=0$ for all $f\in (R^{O(5)}_{-})^{\#}$.
\end{lemma} 
\begin{proof} Let $U,V\in\W_{-}^{\#}$. Similarly to the proof of Corollary~\ref{cor1} we have $\si_6(U+V)=0$ and $\si_8(U+V)=0$. Therefore,
\begin{enumerate} 
\item[$\bullet$] $\si_3(UV)\equiv -\tr(U^3 V^3) - \tr(U^2 V^2 U V) - \tr(U^2 V U V^2)$ in $R^{O(5)}_{-}$,  
 
\item[$\bullet$] $\si_4(UV)\equiv -\si_2(U^2 V^2)+\tr(U^4 V^4) + \tr(U^3 V^3 U V) + \tr(U^3 V U V^3)+ \tr(U^2 V^2 U V U V)+ \tr(U^2 V U V^2 U V)+ \tr(U^2 V U V U V^2)$ in $R^{O(5)}_{-}$. 
\end{enumerate}
Since $\sum_{i=0}^5 U^{5-i} \si_{i}(V)=0$, we obtain that $\tr(U^5 V)\equiv 0$ in $R^{O(5)}_{-}$. By Theorem~\ref{theo_matrix} and Lemma~\ref{lemma_key2}, to complete the proof it remains to note that if $i+j$ is odd, then $\tr(A_1^i A_2^j)=0$.
\end{proof}  

\begin{proof_of}{of part~d) of Theorem~\ref{theo_hsp}.} 
Let $A_1,A_2\in S_{-}$ and $f(A_1,A_2)=0$ for all $f\in H$. By Lemmas~\ref{lemma_key} and~\ref{lemma_n5_1}, to prove the required statement it is enough to show that conditions~a) and~b) from Lemma~\ref{lemma_n5_1} hold.

If $A_1$ or $A_2$ is zero, then the above mentioned conditions are valid. Assume that $A_1\neq0$ and $A_2\neq0$. By Lemma~\ref{lemma_canon}, we can assume that $A_1\in\{Q_1,Q_2,Q_3\}$ and $A_2=\Skew(a_2,b_2,c_2,d_2,e_2,f_2,g_2,h_2,i_2,j_2)$.

Let $A_1=Q_1$. Equalities $\tr(A_1 A_2)=0$ and $\tr(A_1^2 A_2^2)=0$ imply $a_2=-\ii e_2$ and $d_2 = \ii i_2 + \de \ii (c_2 - \ii h_2)$, where $\de\in\{-1,1\}$. 
Since $\tr(A_1^2 A_2^4)=0$, one of the  following possibilities holds.
\begin{enumerate}
\item[$\bullet$] Let $f_2=-\de \ii g_2$. Then $A_1^3=0$ and $A_1 A_2^j A_1=0$ for $1\leq j\leq 4$. 

\item[$\bullet$] Let $c_2=\ii h_2$. Then equalities $\si_2(A_2)=\si_4(A_2)=0$ imply $b_2=0$ and $f_2^2 + g_2^2 + j_2^2 =0$. Hence $A_1 A_2 A_1 = A_2 A_1 A_2 = 0$ and $A_1^2 A_2^2 = 0$.   
\end{enumerate}
In both cases conditions~a) and~b) from Lemma~\ref{lemma_n5_1} hold.

Let $A_1=Q_3$. Then $\tr(A_1^3 A_2)=0$ implies $g_2=\ii a_2 + c_2 +\ii j_2$. Since $\tr(A_1^4 A_2^2)=0$, we have $b_2 = -\ii i_2$. Thus the equality $\tr(A_1 A_2)=0$ implies $j_2 = \ii c_2 + \frac{1+\ii}{2} e_2 - \frac{1-\ii}{2} h_2$. Considering $\tr(A_1^2 A_2^2)=0$, we can see that $e_2=-\ii h_2$. Finally, using $\si_2(A_2)=\si_4(A_2)=0$, we obtain $d_2=f_2=0$. It is not difficult to verify that $\tr(A_1^3 A_2^3)=0$, $A_1 A_2 A_1 A_2 A_1 =0$, and $A_1^{i_1} A_2^{j_1} A_1^{i_2} A_2^{j_2}=0$ for all $i_1,j_1,i_2,j_2\in\{1,\ldots, 4\}$ with $i_1+j_1+i_2+j_2>4$. Thus, conditions~a) and~b) from Lemma~\ref{lemma_n5_1} hold.
   
Let $A_1= Q_2$. If $A_2\sim Q_1$ or $A_2\sim Q_3$, then considering $A_2$ instead of $A_1$ we obtain the required. Assume that $A_2\sim Q_2$. Since $\tr(A_1 A_2)=0$, we have $a_2=-\ii b_2 - \ii f_2 + h_2$. Thus $A_1^2=A_1 A_2 A_1=0$. Similarly we obtain $A_2^2=A_2 A_1 A_2=0$. Therefore, conditions~a) and~b) from Lemma~\ref{lemma_n5_1} hold. The proof is completed.
\end{proof_of}

\section{Polynomial algebras}\label{section_poly}

As an application, we prove the following lemma.

\begin{lemma}\label{lemma_poly}
Let $n\leq 5$. Then the following conditions are equivalent:
\begin{enumerate}
\item[a)] $R_{-}^{O(n)}$ is a polynomial algebra (i.e.~free algebra over $\FF$);

\item[b)] $d=1$ or $(n,d)$ is equal to $(3,2)$. 
\end{enumerate}
\end{lemma}
\begin{proof} By part~a) of Theorem~\ref{theo_hsp}, Theorem~\ref{theo_mgs}, and Lemmas~\ref{lemma_n2},~\ref{lemma_d1}, the required statement holds for $n\leq 3$ as well as for $d=1$.

Consider the case of $n=4$ and $d=2$. Assume that there are $\al_1,\al_2,\be_1,\be_2\in\FF$ such that
\begin{eq}\label{eq_trXXYY_1}
\al_1 \tr(Y_1^2 Y_2^2) + \al_2 \si_2(Y_1 Y_2) 
- \be_1 \si_2(Y_1) \si_2(Y_2) - \be_2 \tr(Y_1 Y_2)^2=0.
\end{eq}%
We set $Y_1=\Skew(a,b,0,0,0,0)$ and $Y_2=\Skew(1,1,0,0,0,0)$ for $a,b\in\FF$ and  consider the left hand side of~\Ref{eq_trXXYY_1} as a polynomial in $a,b$. Hence coefficients of $a^2$ and $ab$ are equal to zero. Thus, $\al_1=\be_1=-2\al_2+4 \be_2$. We set $Y_1=\Skew(a,0,0,0,0,b)$ and $Y_2=\Skew(1,0,0,0,0,1)$ for $a,b\in\FF$. As above, we obtain $\al_1=\al_2=\be_1=\be_2=0$.
This reasoning together with Lemma~\ref{lemma_key2} implies that
\begin{eq}\label{eq_trXXYY_2}
\al_1 \tr(Y_1^2 Y_2^2) + \al_2 \si_2(Y_1 Y_2) \not\equiv0
\end{eq}%
for all $\al_1,\al_2\in\FF$ satisfying $\al_1\neq0$ or $\al_2\neq0$.

Using~\Ref{eq_trXXYY_2}, we can see that every m.h.s.g.~for $R_{-}^{O(4)}$ contain exactly two elements of multidegree $(2,2)$, which are linear combinations of $\tr(Y_1^2 Y_2^2)$ and $\si_2(Y_1 Y_2)$. On the other hand, there is a m.h.s.g.~$S$ for $R_{-}^{O(4)}$ containing the h.s.p.~from part~b) of Theorem~\ref{theo_hsp}. Since the given h.s.p.~contains only one element of multidegree $(2,2)$, we obtain that its cardinality is less than the cardinality of $S$. Thus, $R_{-}^{O(n)}$ is not a polynomial algebra. 

Similarly, we prove this lemma for $n=5$ and $n=2$. 

Note that if $R_{-}^{O(n)}$ is a polynomial algebra for $d=d_1$, then $R_{-}^{O(n)}$ is a polynomial algebra for $d=d_2$, where $d_2<d_1$. This remark completes the proof.  
\end{proof}

\begin{conj} $R_{-}^{O(n)}$ is a polynomial algebra if and only if condition~b) from Lemma~\ref{lemma_poly} holds.
\end{conj}

\end{document}